\documentclass[a4paper]{aomart}

\usepackage{ucs}
\usepackage[utf8x]{inputenc}
\usepackage[english]{babel}

\usepackage[none]{hyphenat}

\usepackage{amsmath}
\usepackage{amsthm}
\usepackage{amsfonts, amssymb}

\input{principio_variacional_para_pressao/definitions}
\DeclareMathOperator{\Ppressureop}{P}
\DeclareMathOperator{\Spressureop}{S}
\DeclareMathOperator{\Gpressureop}{G}
\DeclareMathOperator{\Qpressureop}{Q}

\newcommand{\measurepressure}[3]{{P_{#1}\left({#2}, {#3}\right)}}
\newcommand{\topologicalpressure}[2]{{\Ppressureop\left({#1}, {#2}\right)}}
\newcommand{\dtopologicalpressure}[3]{{\Ppressureop_{#1}\left({#2}, {#3}\right)}}
\newcommand{\separatedpressure}[3]{{\Spressureop_{#1}\left({#2}, {#3}\right)}}
\newcommand{\generatingpressure}[3]{{\Gpressureop_{#1}\left({#2}, {#3}\right)}}
\newcommand{\plustopologicalpressure}[2]{{\Qpressureop^+\left({#1}, {#2}\right)}}
\newcommand{\minustopologicalpressure}[2]{{\Qpressureop^-\left({#1}, {#2}\right)}}

\newcommand{\Pcoverpressure}[4]{{\Ppressureop_{#1}\left({#2}, {#3}, {#4}\right)}}
\newcommand{\Qcoverpressure}[4]{{\Qpressureop_{#1}\left({#2}, {#3}, {#4}\right)}}
\newcommand{\separatedpathpressure}[5][]{{\Spressureop^{#1}_{#2}\left({#3}, {#4}, {#5}\right)}}
\newcommand{\generatingpathpressure}[5][]{{\Gpressureop^{#1}_{#2}\left({#3}, {#4}, {#5}\right)}}

\newcommand{\tpartitionpressure}[4]{{P_{#1}\left({#2}, {#3}, {#4}\right)}}
\newcommand{\tPcoverpressure}[3]{{\Ppressureop\left({#1}, {#2}, {#3}\right)}}
\newcommand{\plustcoverpressure}[3]{{\Qpressureop^+\left({#1}, {#2}, {#3}\right)}}
\newcommand{\minustcoverpressure}[3]{{\Qpressureop^-\left({#1}, {#2}, {#3}\right)}}

\title{Topological Pressure
       \\ for Locally Compact Metrizable Systems}

\author[André Caldas]{André Caldas de Souza}
\email{andrecaldas@unb.br}
\address{Departamento de Matemática -- Universidade de Brasília, Brasília, Brasil}
\fancypagestyle{firstpage}{\fancyhf{}\chead{}}

\begin{document}
  \maketitle

  {

  {

\begin{abstract}
  It is widely known that when
  $X$ is compact Hausdorff,
  and when
  $\dynamicalsystem{X}{T}$ and $\function{f}{X}{\reals}$
  are continuous,
  \begin{equation*}
    \topologicalpressure{T}{f}
    =
    \sup_{\text{$\mu$: Radon probability}}
    \left(
      \measureentropy{\mu}{T}
      +
      \integral{f}{\mu}
    \right),
  \end{equation*}
  where $\topologicalpressure{T}{f}$ is the \emph{topological pressure}
  and $\measureentropy{\mu}{T}$ is the measure theoretic entropy of $T$
  with respect to $\mu$.
  This result is known as \emph{variational principle}.

  We generalize the concept of
  \emph{topological pressure}
  for the case where $X$ is a separable locally compact metric space.
  Our definitions are quite similar to those used in the compact case.
  Our main result is the validity of the
  \emph{variational principle}
  (Theorem \ref{th:principio_variacional}).
\end{abstract}
}
  {

\section{Introduction}

  Traditionally,
  as happens with topological entropy,
  \emph{topological pressure} has been applied for
  dynamical systems defined over compact spaces.
  As a special case,
  for compact metric spaces.
  In this paper we generalize the concept of
  \emph{topological pressure}
  for locally compact subsystems of compact metric ones
  (Definition \ref{def:topological_pressure}).
  Another way to state this condition
  is to say that the system is defined over a
  locally compact separable metric space.
  Or yet,
  that the refered space possesses a one-point metrizable compactification.
  Our definition is quite similar to the compact case,
  and our main result is the validity of the
  \emph{variational principle}
  (Theorem \ref{th:principio_variacional}).
  A consequence of having a one-point metrizable compactification
  is the fact that every Borel measure is in fact a Radon measure
  (see Chapter II, Theorem 3.2 from \cite{probabilities_on_metric_spaces}).

  We try to follow the same approach used for the compact case,
  presented in \cite{viana:ergodic_theory}.
  When $X$ is a locally compact separable metrizable space,
  it possesses a one-point compactification $\onepointcompactification{X}$ with a metric $d$.
  It does not mean that a system $\dynamicalsystem{X}{T}$
  can be continuously extended to a system
  $\function{T^*}{\onepointcompactification{X}}{\onepointcompactification{X}}$.
  However,
  it can be extended to
  $\dynamicalsystem{Z}{S}$,
  where $Z \supset X$ is compact with a metric $r$,
  and $X$ has the topology induced from $Z$
  (Lemma \ref{lemma:compactificacao}).
  In this context, we have:
  \begin{enumerate}
    \item
      The spaces $X$, $\onepointcompactification{X}$ and $Z$.

    \item
      The dynamical systems $T$ and $S$.

    \item
      The metrics $d$, $d$ restricted to $X$,
      $r$ and $r$ restricted to $X$.

    \item
      A weak* sequentially compact space of Borel probabilities over $Z$,
      and its restriction to $X$,
      giving a set of measures $\mu$ such that
      $0 \leq \mu(X) \leq 1$.

    \item
      The spaces of continuous functions over $X$, $\onepointcompactification{X}$ and $Z$.
  \end{enumerate}
  Instead of studying directly the system $T$,
  we shall look at $\dynamicalsystem{Z}{S}$.
  However,
  since we do not want to capture any complexity for $S$ over $Z \setminus X$,
  we shall not make use of the metric $r$.
  When defining \emph{topological pressure},
  instead of $r$,
  we make use of $d$ restricted to $X$.
  Also,
  for the measure theoretic pressure,
  Lemma \ref{lemma:entropia_na_compactificacao}
  allows us to make calculations avoiding the complexity of $S$ outside $X$.
}
  {

\section{Preliminaries}
  \label{section:preliminaries}

  This section is devoted to recalling some elementary definitions
  related to different types of pressure,
  and to proving some fundamental facts which are used in the sequel.
  We also extend the concept of \emph{topological pressure},
  originally defined only for compact systems.

  {

  A \emph{topological dynamical system}
  ---
  or simply a \emph{dynamical system}
  ---
  $\dynamicalsystem{X}{T}$ is a continuous map $T$
  defined over a topological space $X$.
  A \emph{measurable dynamical system}
  $\dynamicalsystem{X}{T}$ is a measurable map $T$
  defined over a measurable space $X$.
  If we embed $X$ with the Borel $\sigma$-álgebra,
  a \emph{topological dynamical system}
  becomes also a \emph{measurable dynamical system}.

  Recall that a family $\family{A}$ of
  subsets of $X$ is a \emph{cover of $X$}
  when
  \begin{equation*}
    X
    =
    \bigcup_{A \in \family{A}} A.
  \end{equation*}
  If the sets in $\family{A}$ are disjoint, then we say that
  $\family{A}$ is a \emph{partition of $X$}.
  A \emph{subcover} of $\family{A}$ is a family
  $\family{B} \subset \family{A}$ which is itself a cover of $X$.
  If $\family{A}$ is a cover of $X$ and $Y \subset X$,
  then we denote by $Y \cap \family{A}$ the cover of $Y$ given by
  \begin{equation*}
    Y \cap \family{A}
    =
    \setsuchthat{A \cap Y}{A \in \family{A}}.
  \end{equation*}

  Given two covers $\family{A}$ and $\family{B}$ of an arbitrary
  set $X$, we say that $\family{A}$
  is \emph{finer} then $\family{B}$ or that
  $\family{A}$ \emph{refines} $\family{B}$
  --- and write $\family{B} \prec \family{A}$ ---
  when every element of $\family{A}$
  is a subset of some element of $\family{B}$.
  We also say that $\family{B}$ is \emph{coarser} then $\family{A}$.
  The relation $\prec$ is a \emph{preorder}, and if we identify the
  \emph{symmetric} covers
  (i.e.: covers $\family{A}$ and $\family{B}$
  such that $\family{A} \prec \family{B}$
  and $\family{B} \prec \family{A}$),
  we have a \emph{lattice}.
  As usual,
  $\family{A} \vee \family{B}$ denotes the representative of the coarsest
  covers of $X$ that refines both $\family{A}$ and $\family{B}$ given by
  \begin{equation*}
    \family{A} \vee \family{B}
    =
    \setsuchthat{A \cap B}{A \in \family{A}, B \in \family{B}, A \cap B \neq \emptyset}.
  \end{equation*}

  Given a dynamical system $\dynamicalsystem{X}{T}$
  and a cover $\family{A}$,
  for each $n \in \naturals$ we define
  \begin{equation*}
    \family{A}^n
    =
    \family{A}
    \vee
    T^{-1}(\family{A})
    \vee \dotsb \vee
    T^{-(n-1)}(\family{A})
  \end{equation*}
  If we want to emphasise the dynamical system $T$,
  we write $\family{A}_T^n$ instead.
  And if $\function{f}{X}{\reals}$ is any function,
  we shall denote by
  \begin{equation*}
    f_n
    =
    \sum_{j=0}^{n-1}
    f \circ T^j.
  \end{equation*}
  And if we really need to emphasise $T$,
  we write $f_{T,n}$ instead.

  For any pseudometric space $\metricspace{X}{d}$,
  given $\varepsilon > 0$ and $x \in X$,
  we denote by
  \begin{equation*}
    \ball{d}{\varepsilon}{x}
    =
    \setsuchthat{y \in X}{d(x,y) < \varepsilon}
  \end{equation*}
  the \emph{open ball} of radius $\varepsilon$,
  centered at $x$.
  And
  \begin{equation*}
    \balls{d}{\varepsilon}
    =
    \setsuchthat{\ball{d}{\varepsilon}{x}}{x \in X}
  \end{equation*}
  is the cover of $X$ composed of all \emph{open balls}
  with radius $\varepsilon$.
}
  {\subsection{Compactification}

  We are mainly interested in dynamical systems
  defined over a metrizable locally compact separable space $X$.
  This is the same as requiring $X$
  to have a metrizable one-point compactification.
  And in general,
  this is not the same as requiring
  that the system $T$ can itself
  be continuously extended to the one-point compactification of $X$.
  The topology of $X$ can be induced by different metrics.
  In special,
  it can be induced by a metric restricted from its one-point compactification.

  To demonstrate our extended version of the \emph{variational principle} for pressure,
  we shall regard the dynamical system
  $\dynamicalsystem{X}{T}$ as a subsystem of a compact metrizable one.
  Please,
  refer to \cite{patrao_caldas:entropia}
  for a detailed treatment
  of the results stated in this subsection.

  \begin{definition}[Subsystem]
    \label{def:subsistema}
    We say that a dynamical system
    $\dynamicalsystem{X}{T}$
    is a \emph{subsystem} of
    $\dynamicalsystem{Z}{S}$
    when $X \subset Z$ has the induced topology and
    $T(x) = S(x)$ for every $x \in X$.
    We also say that
    $S$ \emph{extends} $T$ to $Z$.
  \end{definition}

  \begin{lemma}
    \label{lemma:partition_restriction}
    Suppose $\dynamicalsystem{X}{T}$ is a subsystem of
    $\dynamicalsystem{Z}{S}$.
    If $\family{Z}$ is a covering of $Z$ and
    $\family{C} = X \cap \family{Z}$,
    then
    \begin{equation*}
      \family{C}_T^n
      =
      X
      \cap
      \family{Z}_S^n.
    \end{equation*}
  \end{lemma}

  \begin{proof}
    This is Lemma 2.2 from \cite{patrao_caldas:entropia}.
  \end{proof}

  Under the conditions of
  local compacity, separability and metrizability,
  the dynamical system $\dynamicalsystem{X}{T}$ can always
  be extended to a metrizable compact system.

  \begin{lemma}
    \label{lemma:compactificacao}
    Suppose that $X$ is a topological space
    with metrizable one-point compactification $\onepointcompactification{X} = X \cup \set{\infty}$.
    Then,
    any topological dynamical system $\dynamicalsystem{X}{T}$
    can be extended to a dynamical system $\dynamicalsystem{Z}{S}$,
    with $Z$ compact metrizable,
    and such that the natural projection
    \begin{equation*}
      \functionarray{\pi}{Z}{\onepointcompactification{X}}{x}
                    {\pi(x) = \begin{cases}x, &x \in X\\\infty, &x \not \in X\end{cases}}
    \end{equation*}
    is continuous.
  \end{lemma}

  \begin{proof}
    This is Lemma 2.3 from
    \cite{patrao_caldas:entropia}.
  \end{proof}

  This projection $\pi$,
  on Lemma \ref{lemma:compactificacao},
  induces the pseudometric
  \begin{equation*}
    \widetilde{d}(x,y)
    =
    d(\pi(x), \pi(y))
  \end{equation*}
  over $Z$.
  We denote by the same letter $d$
  the metric over $\onepointcompactification{X}$ and its restriction to $X$.
  Since $\pi$ is continuous,
  this pseudometric $\widetilde{d}$ is such that the ``open balls''
  are in fact open in the topology of $Z$,
  although in general,
  $\widetilde{d}$ does not generate the topology.
  Denote by
  \begin{equation*}
    \complementset{X}
    =
    Z \setminus X
  \end{equation*}
  the complement of $X$ in $Z$.
  Since the set $\complementset{X}$ has zero diameter with respect to $\widetilde{d}$,
  the balls
  $\ball{\widetilde{d}}{\varepsilon}{z}$
  either contain $\complementset{X}$,
  or have empty intersection with it.

  \begin{definition}[One-Point Metric]
    Whenever $X$ has a metrizable one-point compactification $\onepointcompactification{X}$,
    we shall call the restriction of a metric $d$ over $\onepointcompactification{X}$ to $X$
    a \emph{one-point metric}.
  \end{definition}

  In particular,
  if we say that $X$ has a \emph{one-point metric} $d$,
  this implies that $X$ has a one-point compactification.

  Under the conditions of
  Lemma \ref{lemma:compactificacao},
  $X$ is an open subset of $Z$.
  In fact,
  $\complementset{X} = \pi^{-1}(\infty)$
  is closed.
  In this case,
  the Borel sets of $X$ are Borel sets of $Z$,
  and we may restrict Borel measures over $Z$
  to the Borel sets of $X$
  and produce a Borel measure over $X$.
  On the other hand,
  if $\mu$ is a Borel measure over $X$,
  we can extend it to $Z$ by declaring $\mu(\complementset{X}) = 0$.
  We shall use the same letter $\mu$ to denote a measure over $Z$
  as well as its restriction to $X$.
  If we want to make the distinction clear,
  we may write $\mu|_X$ instead.

  From now on,
  unless explicit mention to the contrary,
  $\dynamicalsystem{X}{T}$ will be a topological dynamical system,
  where $X$ admits a one-point metrizable compactification $\onepointcompactification{X} = X \cup \set{\infty}$.
  Also,
  $\dynamicalsystem{Z}{S}$ will be a continuous extension of $T$
  and $\function{\pi}{Z}{\onepointcompactification{X}}$ the natural projection
  whose properties and existence are asured by
  Lemma \ref{lemma:compactificacao}.
  Also,
  $d$ will be a one-point metric and
  $\widetilde{d}(x,y) = d(\pi(x), \pi(y))$
  will be the pseudometric induced in $Z$ by $d$.

  \begin{definition}[One-Point Uniformly Continuous]
    \label{definition:one_point_continuous}
    If $\function{f}{X}{\reals}$ is uniformly continuous
    with respect to some one-point metric
    (and therefore, every one-point metric),
    we shall say that $f$ is \emph{one-point uniformly continuous}.
  \end{definition}

  A \emph{one-point uniformly continuous}
  $\function{f}{X}{\reals}$
  is nothing more then the restriction to $X$
  of a continuous $\function{f}{\onepointcompactification{X}}{\reals}$,
  which we shall denote by the same letter $f$.
  We can always write a
  \emph{one-point uniformly continuous} function
  as a sum $f + c$,
  where $f$ \emph{vanishes at infinity}
  (i.e.: $f \in \vanishingatinfinityfunctions{X}$)
  and $c \in \reals$ is a constant.

  The \emph{one-point uniformly continuous} $f$
  induces the continuous $g = f \circ \pi$,
  from $Z$ to $\reals$.
  Notice that
  since we have the dynamical systems
  $\dynamicalsystem{X}{T}$
  and
  $\dynamicalsystem{Z}{S}$,
  and the functions
  $\function{f_n}{X}{\reals}$
  and
  $\function{g_n}{Z}{\reals}$
  are defined.
  That is,
  \begin{align*}
    f_n
    &=
    f + f \circ T + \dotsb + f \circ T^{n-1}
    \\
    g_n
    &=
    g + g \circ S + \dotsb + f \circ S^{n-1}.
  \end{align*}
  However,
  $\function{f_n}{\onepointcompactification{X}}{\reals}$
  is, in principle, not defined.
}
  {\subsection{Pressure with a Measure}

  Now,
  we shall define the concept of \emph{pressure} of a measurable
  dynamical system with respect to an invariant finite measure.
  Traditionally,
  \emph{pressure} has been defined only for probability measures.
  However,
  as shown in \cite{patrao_caldas:entropia},
  extending this to finite measures is straight forward,
  and can be quite useful in topological dynamical systems
  where $X$ is not compact.
  First,
  we recall some definitions.
  More details can be found in \cite{patrao_caldas:entropia}.

  \begin{definition}[Kolmogorov-Sinai Entropy]
    Consider the finite measure space $\measurespace{X}$ and a
    finite measurable partition $\family{C}$.
    The \emph{partition entropy} of $\family{C}$ is
    \begin{equation*}
      \partitionentropy{\mu}{\family{C}}
      =
      \sum_{C \in \family{C}}
      \mu(C) \log \frac{1}{\mu(C)}.
    \end{equation*}
    For the measurable dynamical system $\dynamicalsystem{X}{T}$,
    if $\mu$ is a $T$-invariant finite measure,
    the \emph{partition entropy of $T$ with respect to $\family{C}$} is
    \begin{equation*}
      \tpartitionentropy{\mu}{\family{C}}{T}
      =
      \lim_{n \rightarrow \infty}
      \frac{1}{n} \partitionentropy{\mu}{\family{C}^n},
    \end{equation*}
    and the \emph{Kolmogorov-Sinai entropy} of $T$ is
    \begin{equation*}
      \measureentropy{\mu}{T}
      =
      \sup_{\substack{\family{C} \text{: finite} \\ \text{measurable partition}}}
      \tpartitionentropy{\mu}{\family{C}}{T}.
    \end{equation*}
  \end{definition}

  \begin{definition}[Pressure]
    \label{definition:pressure}
    Let $\dynamicalsystem{X}{T}$ be a measurable dynamical system
    with a finite $T$-invariant measure $\mu$.
    Let $\function{f}{X}{\reals}$ be an integrable function.
    The quantity
    \begin{equation*}
      \measurepressure{\mu}{T}{f}
      =
      \measureentropy{\mu}{T}
      +
      \integral{f}{\mu}
    \end{equation*}
    is the \emph{pressure} of $T$
    with respect to the measure $\mu$
    and potential $f$.
    To make the notation cleaner,
    also define
    \begin{equation*}
      \tpartitionpressure{\mu}{T}{f}{\family{C}}
      =
      \tpartitionentropy{\mu}{\family{C}}{T}
      +
      \integral{f}{\mu}.
    \end{equation*}
  \end{definition}

  {\subsubsection{Properties}

  In our way into showing the \emph{variational principle},
  we shall make use of some techiniques that are already quite standard,
  and some that are not.
  The following fact is not a standard result.
  It was presented in \cite{patrao_caldas:entropia},
  and relates the entropy of a system
  and the entropy of it's extension.

  \begin{lemma}
    \label{lemma:entropia_na_compactificacao}
    Let
    $\dynamicalsystem{Z}{S}$
    be a measurable dynamical system and
    $\dynamicalsystem{X}{T}$
    a subsystem with $X \subset Z$ measurable.
    If $\mu$ is an $S$-invariant measure,
    and if
    \begin{equation*}
      \family{Z}
      =
      \set{Z_0, \dotsc, Z_k}
    \end{equation*}
    is a measurable partition of $Z$
    such that
    $\complementset{X} \subset Z_0$,
    then
    $\mu$ is $T$-invariant and
    \begin{equation*}
      \tpartitionentropy{\mu}{\family{Z}}{S}
      \leq
      \measureentropy{\mu}{T}.
    \end{equation*}
    If
    $f \in \vanishingatinfinityfunctions{X}$,
    then
    \begin{equation*}
      \tpartitionpressure{\mu}{S}{f \circ \pi}{\family{Z}}
      \leq
      \measurepressure{\mu}{T}{f}.
    \end{equation*}
  \end{lemma}

  \begin{proof}
    The first part is Lemma 2.10 from \cite{patrao_caldas:entropia}.
    The second part follows from the fact that
    \begin{equation*}
      \integral{f \circ \pi}{\mu}
      =
      \integral{f}{\mu|_X}
      +
      \mu\left(\complementset{X}\right)
      f(\infty)
      =
      \integral{f}{\mu|_X}.
    \end{equation*}
  \end{proof}

  Since we do not require the $T$-invariant measure to be a probability,
  the following lemma can be quite handy.

  \begin{lemma}
    \label{le:ksentropy_com_medidas_neq_1}
    Given a measurable dynamical system $\dynamicalsystem{X}{T}$ and a finite
    $T$-invariant measure $\mu$,
    then, for $\alpha \geq 0$,
    \begin{equation*}
      \measureentropy{\alpha \mu}{T}
      =
      \alpha \measureentropy{\mu}{T}.
    \end{equation*}
  \end{lemma}

  \begin{proof}
    This is Lemma 2.9 from \cite{patrao_caldas:entropia}.
  \end{proof}

  Now,
  we list some properties of the pressure,
  most of them are just a consequence of some corresponding property
  of the Kolmogorov-Sinai entropy.

  \begin{proposition}
    \label{proposition:sistema_iterado:medida}
    Let $\dynamicalsystem{X}{T}$ be a measurable dynamical system,
    $\mu$ a $T$-invariant finite measure,
    and $\function{f}{X}{\reals}$ an integrable function,
    Then,
    \begin{equation*}
      \measurepressure{\mu}{T^k}{f_k}
      =
      k
      \measurepressure{\mu}{T}{f}.
    \end{equation*}
  \end{proposition}

  \begin{proof}
    Notice that since $\mu$ is $T$-invariant,
    \begin{equation*}
      \integral{f_k}{\mu}
      =
      k
      \integral{f}{\mu}.
    \end{equation*}
    Also,
    The equality
    \begin{equation*}
      \measureentropy{\mu}{T^k}
      =
      k
      \measureentropy{\mu}{T}
    \end{equation*}
    is usualy stated for the case where $\mu$ is a probability measure
    (for example, Proposition 9.1.14 from \cite{viana:ergodic_theory}).
    For the general case,
    where $\mu$ is a finite measure,
    combine this with Lemma \ref{le:ksentropy_com_medidas_neq_1}.
    Of course,
    the original demonstration for probability measure works
    verbatim for the more general finite measure case
    (see Remark 2.20 in \cite{patrao_caldas:entropia}).
  \end{proof}

  When defining \emph{topological pressure}
  we shall use the concept of \emph{admissible cover}
  (Definition \ref{def:cobertura_admissivel}).
  A measure theoretic counterpart is the \emph{admissible partition}.

  \begin{definition}[Admissible Partition]
    In a topological space $\topologicalspace{X}$,
    a finite (measurable) partition is said to be \emph{admissible}
    when every element but one is compact.
  \end{definition}

  The pressure can be calculated using \emph{admissible partitions}.
  This fact is usualy not explicitly stated
  as we did in Proposition \ref{proposition:pressure_with_admissible_partitions}.
  But it is not new, as it is usually embedded in the demonstrations
  of the \emph{variational principle} for the compact case
  (see, for example, \cite{viana:ergodic_theory,walters,patrao_caldas:entropia}).

  \begin{proposition}
    \label{proposition:pressure_with_admissible_partitions}
    If $\dynamicalsystem{X}{T}$ is a topological dynamical system,
    $\mu$ is a $T$-invariant Radon probability measure,
    and $\function{f}{X}{\reals}$ is an integrable function.
    Then,
    \begin{equation*}
      \measurepressure{\mu}{T}{f}
      =
      \sup_{\substack{\family{K} \text{: admissible} \\ \text{partition}}}
      \tpartitionpressure{\mu}{T}{f}{\family{K}}.
    \end{equation*}
  \end{proposition}

  For the proof of Proposition \ref{proposition:pressure_with_admissible_partitions},
  we need the concept of \emph{conditional entropy}.
  The proof will be presented after some preparation.

  \begin{definition}[Conditional Entropy]
    Given a probability measure $\mu$
    and two finite measurable partitions
    $\family{C}$ and $\family{D}$,
    the \emph{conditional entropy}
    is defined as the expected value
    \begin{equation*}
      \conditionalpartitionentropy{\mu}{\family{C}}{\family{D}}
      =
      \sum_{D \in \family{D}}
      \mu(D)
      \partitionentropy{\conditionalprobability{\mu}{D}{\cdot}}{\family{C}}.
    \end{equation*}
  \end{definition}

  Conditional entropy possesses the following properties.

  \begin{lemma}
    \label{lemma:estimation_using_conditional_entropy}
    Let $\dynamicalsystem{X}{T}$ be a measurable dynamical system
    with $T$-invariant probability measure $\mu$.
    If $\family{C}$ and $\family{D}$ are two measurable finite partitions,
    then
    \begin{equation*}
      \tpartitionentropy{\mu}{\family{C}}{T}
      \leq
      \tpartitionentropy{\mu}{\family{D}}{T}
      +
      \conditionalpartitionentropy{\mu}{\family{C}}{\family{D}}.
    \end{equation*}
  \end{lemma}

  \begin{proof}
    This is item $(iv)$ of Theorem 4.12 from \cite{walters}.
    Or Lemma 9.1.11 from \cite{viana:ergodic_theory}.
  \end{proof}

  We shall need to calculate the \emph{conditional entropy}
  only for the following case.

  \begin{lemma}
    \label{lemma:conditional_entropy:from_inside}
    Let $\family{C} = \set{C_1, \dotsc, C_n}$ be a measurable partition.
    If $\family{K} = \set{K_0, K_1, \dotsc, K_n}$ is such that
    $K_j \subset C_j$ for every $j = 1, \dotsc, n$,
    then
    \begin{align*}
      \conditionalpartitionentropy{\mu}{\family{C}}{\family{K}}
      &=
      \mu(K_0)
      \partitionentropy{\conditionalprobability{\mu}{K_0}{\cdot}}{\family{C}}
      \\
      &\leq
      \mu(K_0)
      \log n.
    \end{align*}
  \end{lemma}

  \begin{proof}
    One just has to notice that for every $C \in \family{C}$
    and $j \neq 0$,
    $\conditionalprobability{\mu}{K_j}{C}$
    is either $0$ or $1$.
    Therefore,
    for $j = 1, \dotsc, n$,
    \begin{equation*}
      \partitionentropy{\conditionalprobability{\mu}{K_j}{\cdot}}{\family{C}}
      =
      0.
    \end{equation*}
    It is a very well known fact that
    $\partitionentropy{\nu}{\family{C}} \leq \log n$
    (see, for example, Lemma 9.1.3 from \cite{viana:ergodic_theory}).
  \end{proof}

  We are now,
  ready do demonstrate
  Proposition \ref{proposition:pressure_with_admissible_partitions}.

  \begin{proof}[Proof (Proposition \ref{proposition:pressure_with_admissible_partitions})]
    If $\family{K}$ is an admissible partition,
    it is finite by definition,
    and measurable because compact sets are measurable.
    From the definition of $\measureentropy{\mu}{T}$,
    it is evident that
    \begin{align*}
      \sup_{\substack{\family{K} \text{: admissible} \\ \text{partition}}}
      \tpartitionentropy{\mu}{T}{\family{K}}
      &\leq
      \sup_{\substack{\family{C} \text{: finite} \\ \text{measurable partition}}}
      \tpartitionentropy{\mu}{T}{\family{C}}
      \\
      &=
      \measureentropy{\mu}{T}.
    \end{align*}

    To finish the demonstration,
    we just have to find for any $\varepsilon > 0$
    and any measurable finite partition
    $\family{C} = \set{C_1, \dotsc, C_n}$,
    an admissible partition $\family{K}$ such that
    \begin{equation*}
      \tpartitionentropy{\mu}{T}{\family{C}}
      \leq
      \tpartitionentropy{\mu}{T}{\family{K}}
      +
      \varepsilon.
    \end{equation*}
    To that end,
    let's choose the partition
    $\family{K} = \set{K_0, \dotsc, K_n}$,
    where
    $K_j \subset C_j$ for $j = 1, \dotsc, n$,
    and $\mu(K_0) \leq \frac{\varepsilon}{\log n}$.
    For example,
    since $\mu$ is Radon,
    just choose a compact $K_j \subset C_j$ for each $j = 1, \dotsc, n$,
    such that
    \begin{equation*}
      \mu(C_j \setminus K_j)
      \leq
      \frac{\varepsilon}{n \log n}.
    \end{equation*}
    Since $K_0 = \complementset{(K_1 \cup \dotsb \cup K_n)}$,
    \begin{equation*}
      \mu(K_0)
      =
      \sum_{j=1}^n \mu(C_j \setminus K_j)
      \leq
      \frac{\varepsilon}{\log n}.
    \end{equation*}
    Now,
    using Lemmas
    \ref{lemma:estimation_using_conditional_entropy}
    and
    \ref{lemma:conditional_entropy:from_inside},
    \begin{align*}
      \tpartitionentropy{\mu}{T}{\family{C}}
      &\leq
      \tpartitionentropy{\mu}{T}{\family{K}}
      +
      \frac{\varepsilon}{\log n} \log n
      \\
      &=
      \tpartitionentropy{\mu}{T}{\family{K}}
      +
      \varepsilon.
    \end{align*}
  \end{proof}

  Next,
  we present an upper bound for calculating the \emph{pressure}
  that also motivates the definition of \emph{topological pressure}.
  First,
  notice that
  \begin{equation*}
    \measurepressure{\mu}{T}{f}
    =
    \sup_{\family{C}}
    \lim_{n \rightarrow \infty}
    \left(
      \integral{f}{\mu}
      +
      \frac{1}{n}
      \partitionentropy{\mu}{\family{C}^n}
    \right),
  \end{equation*}
  where the supremum is taken over every measurable finite partition $\family{C}$.

  \begin{lemma}
    \label{lemma:measure_pressure:upper_bound}
    Let $\dynamicalsystem{X}{T}$ be a measurable dynamical system,
    $\mu$ a $T$-invariant probability measure,
    and $\function{f}{X}{\reals}$ an integrable function.
    Then,
    for every finite measurable partition $\family{C}$,
    \begin{equation*}
      \integral{f}{\mu}
      +
      \frac{1}{n}
      \partitionentropy{\mu}{\family{C}^n}
      \leq
      \frac{1}{n}
      \log \sum_{C \in \family{C}^n} \sup \expbase^{f_n(C)}.
    \end{equation*}
  \end{lemma}

  \begin{proof}
    Notice that,
    from the $T$-invariance of $\mu$,
    \begin{equation*}
      \integral{f}{\mu}
      =
      \frac{1}{n}
      \integral{f_n}{\mu}.
    \end{equation*}
    Therefore,
    \begin{align*}
      \integral{f}{\mu}
      +
      \frac{1}{n}
      \partitionentropy{\mu}{\family{C}^n}
      &=
      \frac{1}{n}
      \left(
        \integral{f_n}{\mu}
        +
        \sum_{C \in \family{C}^n}
        \mu(C) \log \frac{1}{\mu(C)}
      \right)
      \\
      &\leq
      \frac{1}{n}
      \sum_{C \in \family{C}^n}
      \left(
        \mu(C) \sup f_n(C)
        +
        \mu(C) \log \frac{1}{\mu(C)}
      \right).
    \end{align*}
    Now,
    the result follows from Lemma 10.4.4 from \cite{viana:ergodic_theory}.
  \end{proof}
}
}
  {\subsection{Topological Pressure}
  \label{section:topological_pressure}

  As it happens with \emph{topological entropy} in the compact case,
  there are different equivalent ways to define \emph{topological pressure}.
  For non-compact systems,
  those different definitions might not be equivalent.
  In the same spirit of that from \cite{patrao_caldas:entropia},
  we shall adapt some of those definitions so they
  work in the non-compact case as well.
  As for the notation,
  we try to follow as closely as possible that of
  \cite{viana:ergodic_theory}.
  As in \cite{patrao_caldas:entropia},
  we use \emph{admissible covers} and \emph{covers of balls}
  in order to define the different concepts of \emph{topological pressure}.

  \begin{definition}
    \label{def:potential_cardinality}
    Given
    $\function{f}{X}{\reals}$
    and
    a cover $\family{A}$ of a set $X$,
    define
    \begin{align*}
      \Qcoverpressure{n}{T}{f}{\family{A}}
      &=
      \inf
      \setsuchthat{\sum_{A \in \family{A}'} \inf \expbase^{f_n(A)}}
                  {\text{$\family{A}'$ is a subcover of $\family{A}^n$}}
      \\
      \Pcoverpressure{n}{T}{f}{\family{A}}
      &=
      \inf
      \setsuchthat{\sum_{A \in \family{A}'} \sup \expbase^{f_n(A)}}
                  {\text{$\family{A}'$ is a subcover of $\family{A}^n$}}.
    \end{align*}
  \end{definition}

  The role played by
  $\Qcoverpressure{n}{T}{f}{\family{A}}$
  and
  $\Pcoverpressure{n}{T}{f}{\family{A}}$
  in Definition \ref{def:potential_cardinality}
  is analogous to that of $\covercardinality{\family{A}^n}$
  when we define topological entropy.
  In fact,
  \begin{equation*}
    \Qcoverpressure{n}{T}{0}{\family{A}}
    =
    \Pcoverpressure{n}{T}{0}{\family{A}}
    =
    \covercardinality{\family{A}^n}.
  \end{equation*}
  The following Lemma shows that
  $\Qcoverpressure{n}{T}{f}{\family{A}}$
  has a property very simmilar to that of
  $\covercardinality{\family{A}^n}$.

  \begin{lemma}
    \label{lemma:refinement:inequality}
    If
    $\family{A} \prec \family{B}$,
    then, for any
    $\function{f}{X}{\reals}$,
    and any
    $n = 1, 2, \dotsc$,
    \begin{equation*}
      \Qcoverpressure{n}{T}{f}{\family{A}}
      \leq
      \Qcoverpressure{n}{T}{f}{\family{B}}.
    \end{equation*}
  \end{lemma}

  \begin{proof}
    Notice that
    $\family{A} \prec \family{B}$
    implies
    $\family{A}^n \prec \family{B}^n$.

    For every $B \in \family{B}^n$,
    there is an $A_B \in \family{A}^n$
    such that $B \subset A_B$.
    In this case,
    \begin{equation*}
      \inf \expbase^{f(A_B)}
      \leq
      \inf \expbase^{f(B)}.
    \end{equation*}
    Notice that
    for every subcover $\family{B}'$ of $\family{B}^n$,
    \begin{equation*}
      \family{A}'
      =
      \setsuchthat{A_B \in \family{A}^n}{B \in \family{B}'}
    \end{equation*}
    is a subcover of $\family{A}^n$.
    Therefore,
    \begin{align*}
      \Qcoverpressure{n}{T}{f}{\family{A}}
      &\leq
      \sum_{A \in \family{A}'}
      \inf \expbase^{f(A)}
      \\
      &\leq
      \sum_{B \in \family{B}'}
      \inf \expbase^{f(A_B)}
      \\
      &\leq
      \sum_{B \in \family{B}'}
      \inf \expbase^{f(B)}.
    \end{align*}
    The result follows
    if we take the infimum over every subcover $\family{B}'$ of $\family{B}^n$.
  \end{proof}

  While $\Qcoverpressure{n}{T}{f}{\family{A}}$
  has the property stated in Lemma \ref{lemma:refinement:inequality},
  the sequence $\Pcoverpressure{n}{T}{f}{\family{A}}$
  shares a different property with $\covercardinality{\family{A}^n}$:
  it is submultiplicative.
  That is,
  \begin{equation*}
    \Pcoverpressure{m+n}{T}{f}{\family{A}}
    \leq
    \Pcoverpressure{m}{T}{f}{\family{A}}
    \Pcoverpressure{n}{T}{f}{\family{A}}.
  \end{equation*}
  And therefore,
  \begin{equation*}
    \lim_{n \rightarrow \infty}
    \frac{1}{n}
    \log \Pcoverpressure{n}{T}{f}{\family{A}}
  \end{equation*}
  exists.

  \begin{lemma}
    For any cover $\family{A}$ of a set $X$,
    \begin{equation*}
      \lim_{n \rightarrow \infty}
      \frac{1}{n}
      \log \Pcoverpressure{n}{T}{f}{\family{A}}
    \end{equation*}
    exists.
  \end{lemma}

  \begin{proof}
    This is Lemma 9.3 from \cite{walters}.
    But it is important to notice that
    although Walters assumes $X$ to be compact,
    Lemma 9.3 does not depend on this hypothesis.
    It is also worth noticing that
    the demonstration also does not depend on the fact that
    $\family{A}$ is an open cover,
    or that $\function{f}{X}{\reals}$ is continuous,
    and these hypothesis could be removed from the statement of Lemma 9.3.
  \end{proof}

  \begin{definition}
    \label{def:cover_pressure}
    Given
    $\function{f}{X}{\reals}$
    and
    a cover $\family{A}$ of a set $X$,
    define
    \begin{align*}
      \minustcoverpressure{T}{f}{\family{A}}
      &=
      \liminf_{n \rightarrow \infty}
      \frac{1}{n}
      \log \Qcoverpressure{n}{T}{f}{\family{A}}
      \\
      \plustcoverpressure{T}{f}{\family{A}}
      &=
      \limsup_{n \rightarrow \infty}
      \frac{1}{n}
      \log \Qcoverpressure{n}{T}{f}{\family{A}}
      \\
      \tPcoverpressure{T}{f}{\family{A}}
      &=
      \lim_{n \rightarrow \infty}
      \frac{1}{n}
      \log \Pcoverpressure{n}{T}{f}{\family{A}}.
    \end{align*}
  \end{definition}

  As in \cite{patrao_caldas:entropia},
  we shall restrict our attention to
  \emph{admissible covers}.

  \begin{definition}[Admissible Cover]
    \label{def:cobertura_admissivel}
    In a topological space $\topologicalspace{X}$,
    an open cover $\family{A}$ is said to be \emph{admissible}
    when at least one of its elements has compact complement.
    If every set has compact complement,
    $\family{A}$ is said to be \emph{strongly admissible},
    or \emph{s-admissible} for short.
  \end{definition}

  \begin{lemma}
    \label{lemma:cover_from_partition}
    In a topological space $\topologicalspace{X}$,
    if
    \begin{equation*}
      \family{K}
      =
      \set{K_0, \dotsc, K_n}
    \end{equation*}
    is an admissible partition where $K_1, \dotsc, K_n$ are all compact,
    then
    \begin{equation*}
      \family{A}
      =
      \set{K_0 \cup K_1, K_0 \cup K_2, \dotsc, K_0 \cup K_n}
    \end{equation*}
    is a strongly admissible cover.
  \end{lemma}

  \begin{proof}
    One just has to notice that $\family{A}$ does cover $X$.
    And also, that
    \begin{equation*}
      \complementset{(K_0 \cup K_j)}
      =
      \bigcup_{i \in \set{1, \dotsc, n} \setminus \set{j}}
      K_i
    \end{equation*}
    is compact for every $j = 1, \dotsc, n$.
  \end{proof}

  The following Lemma is usually embedded in the demonstration of the \emph{variational principle}.
  It is not new,
  except for the fact that it is usually applied without being formally stated.

  \begin{lemma}
    \label{lemma:cover_from_partition:intersection_counting}
    In a topological space $\topologicalspace{X}$,
    let
    \begin{equation*}
      \family{K}
      =
      \set{K_0, \dotsc, K_n}
    \end{equation*}
    be an admissible partition where $K_1, \dotsc, K_n$ are all compact.
    Let
    \begin{equation*}
      \family{A}
      =
      \set{K_0 \cup K_1, K_0 \cup K_2, \dotsc, K_0 \cup K_n}.
    \end{equation*}
    If $\family{B}$ refines $\family{A}$,
    then,
    for each $B \in \family{B}^k$,
    the number of elements of $\family{K}^k$
    that $B$ intersects is at most $2^k$.
  \end{lemma}

  \begin{proof}
    Since $\family{A} \prec \family{B}$,
    $\family{A}^k \prec \family{B}^k$.
    Therefore,
    $B$ is contained in some $A \in \family{A}^k$.
    Now,
    $A$ is of the form
    \begin{equation*}
      (K_0 \cup K_{\lambda_1})
      \cap
      T^{-1}(K_0 \cup K_{\lambda_2})
      \cap \dotsb \cap
      T^{-(n-1)}(K_0 \cup K_{\lambda_k}),
    \end{equation*}
    for some $\lambda \in \set{1, \dotsc, n}^k$.

    Therefore,
    \begin{equation*}
      B
      \subset
      \bigcup_{\gamma \in \set{0,1}^k}
      \left(
        K_{\gamma_1 \lambda_1}
        \cap
        T^{-1} K_{\gamma_2 \lambda_2}
        \cap \dotsb \cap
        T^{-(n-1)} K_{\gamma_n \lambda_n}
      \right).
    \end{equation*}
    Since $\family{K}^k$ partitions $X$,
    $B$ intersects only the non empty sets in this union.
    And since there is one for each $\gamma \in \set{0,1}^k$,
    the claim follows.
  \end{proof}

  An important feature of \emph{admissible covers}
  is that there is a \emph{Lebesgue Number} associated to them.

  \begin{lemma}[Lebesgue Number]
    \label{lemma:lebesgue_number}
    Let $d$ be the restirction to $X$ of a metric in some compactification,
    and let $\family{A}$ be an admissible cover.
    Then, there exists $\varepsilon > 0$ such that
    \begin{equation*}
      \family{A}
      \prec
      \balls{d}{\varepsilon}.
    \end{equation*}
  \end{lemma}

  \begin{proof}
    Remark 2.15
    and
    Lemma 2.27,
    both from \cite{patrao_caldas:entropia},
    lead to the desired result.
  \end{proof}

  \begin{definition}[Topological Pressures]
    For a dynamical system $\dynamicalsystem{X}{T}$,
    and a function $\function{f}{X}{\reals}$,
    define
    \begin{align*}
      \minustopologicalpressure{T}{f}
      &=
      \sup_{\text{$\family{A}$: admissible cover}}
      \minustcoverpressure{T}{f}{\family{A}}
      \\
      \plustopologicalpressure{T}{f}
      &=
      \sup_{\text{$\family{A}$: admissible cover}}
      \plustcoverpressure{T}{f}{\family{A}}.
    \end{align*}
    And if $d$ is a metric over $X$,
    define
    \begin{equation*}
      \dtopologicalpressure{d}{T}{f}
      =
      \limsup_{\varepsilon \rightarrow 0}
      \tPcoverpressure{T}{f}{\balls{d}{\varepsilon}}.
    \end{equation*}
  \end{definition}

  \begin{lemma}
    \label{lemma:balls:one_point_admissible}
    \label{lemma:topological_and_metric:hard_inequality}
    If $d$ is a one-point metric for $X$,
    then
    $\balls{d}{\varepsilon}$ is \emph{admissible} for any $\varepsilon > 0$.
    Also,
    for any $\function{f}{X}{\reals}$,
    \begin{align*}
      \minustopologicalpressure{T}{f}
      &=
      \sup_{\varepsilon > 0}
      \minustcoverpressure{T}{f}{\balls{d}{\varepsilon}}
      \\
      \plustopologicalpressure{T}{f}
      &=
      \sup_{\varepsilon > 0}
      \plustcoverpressure{T}{f}{\balls{d}{\varepsilon}}.
    \end{align*}
    And if $f$ is uniformly continuous with respect to $d$,
    \begin{equation*}
      \dtopologicalpressure{d}{T}{f}
      \leq
      \minustopologicalpressure{T}{f}.
    \end{equation*}
  \end{lemma}

  \begin{proof}
    Just take any $x \in X$ such that
    $d(x, \infty) < \varepsilon$.
    Then,
    \begin{equation*}
      X
      \setminus
      \ball[d]{x}{\varepsilon}
      =
      \onepointcompactification{X}
      \setminus
      \ball[d]{x}{\varepsilon}
    \end{equation*}
    is closed in $\onepointcompactification{X}$,
    and therefore, compact.
    Therefore,
    $\balls{d}{\varepsilon}$ is admissible.

    In particular,
    the definition of
    $\minustopologicalpressure{T}{f}$
    and
    $\plustopologicalpressure{T}{f}$
    implies that
    \begin{align*}
      \sup_{\varepsilon > 0}
      \minustcoverpressure{T}{f}{\balls{d}{\varepsilon}}
      &\leq
      \minustopologicalpressure{T}{f}
      \\
      \sup_{\varepsilon > 0}
      \plustcoverpressure{T}{f}{\balls{d}{\varepsilon}}
      &\leq
      \plustopologicalpressure{T}{f}.
    \end{align*}

    On the other hand,
    if $\family{A}$ is admissible,
    Lemma \ref{lemma:lebesgue_number}
    gives $\varepsilon > 0$ such that
    $\family{A} \prec \balls{d}{\varepsilon}$.
    Therefore,
    Lemma \ref{lemma:refinement:inequality} implies that
    \begin{align*}
      \minustcoverpressure{T}{f}{\family{A}}
      &\leq
      \sup_{\varepsilon > 0}
      \minustcoverpressure{T}{f}{\balls{d}{\varepsilon}}
      \\
      \plustcoverpressure{T}{f}{\family{A}}
      &\leq
      \sup_{\varepsilon > 0}
      \plustcoverpressure{T}{f}{\balls{d}{\varepsilon}}.
    \end{align*}
    By taking the supremum over all admissible covers $\family{A}$,
    \begin{align*}
      \minustopologicalpressure{T}{f}
      &\leq
      \sup_{\varepsilon > 0}
      \minustcoverpressure{T}{f}{\balls{d}{\varepsilon}}
      \\
      \plustopologicalpressure{T}{f}
      &\leq
      \sup_{\varepsilon > 0}
      \plustcoverpressure{T}{f}{\balls{d}{\varepsilon}}.
    \end{align*}

    Finally,
    if $f$ is uniformly continuous with respect to $d$,
    then,
    for each $\eta > 0$,
    there exists $\varepsilon_0 > 0$ such that
    for every $n = 1, 2, \dotsc$,
    every non null $\varepsilon \leq \varepsilon_0$,
    and every $B \in \balls{d}{\varepsilon}^n$,
    \begin{equation*}
      \sup f_n(B)
      \leq
      \inf f_n(B) + n \eta.
    \end{equation*}
    One just has to choose $\varepsilon_0 > 0$ such that
    \begin{equation*}
      d(x,y) < 2 \varepsilon_0
      \Rightarrow
      \abs{f(x) - f(y)} < \eta.
    \end{equation*}
    In this case,
    for any subcover $\family{B} \subset \balls{d}{\varepsilon}^n$
    \begin{equation*}
      \sum_{B \in \family{B}}
      \sup \expbase^{f_n(B)}
      \leq
      \expbase^{n \eta}
      \sum_{B \in \family{B}}
      \inf \expbase^{f_n(B)}.
    \end{equation*}
    Taking the infimum for every subcover $\family{B}$,
    taking the logarithm,
    dividing by $n$
    and taking the $\liminf$,
    \begin{equation*}
      \dtopologicalpressure{d}{T}{f}
      \leq
      \minustopologicalpressure{T}{f}
      +
      \eta.
    \end{equation*}
    Since $\eta$ is arbitrary,
    the result follows.
  \end{proof}

  As in the case of topological entropy,
  we can define yet another concept of topological pressure
  using $(n,\varepsilon)$-separated
  and $(n,\varepsilon)$-generating sets.
  In the compact case,
  those concepts are all equivalent to the ones we have already defined.
  Given a metric $d$ over $X$ and $\varepsilon > 0$,
  we say that a set $E_n$ is $(n, \varepsilon)$-separated
  if for any $x, y \in E_n$,
  \begin{equation*}
    \forall j = 0, \dotsc, n-1,\;
    d(T^j x, T^j y) < \varepsilon
    \Rightarrow
    x = y.
  \end{equation*}
  And we say that $G_n$ is $(n, \varepsilon)$-generating
  if given any $x \in X$,
  there is $y \in G_n$ such that
  for any $j = 0, \dotsc, n-1$,
  \begin{equation*}
    d(T^j x, T^j y) < \varepsilon.
  \end{equation*}
  More information about the relation between
  $(n, \varepsilon)$-separated sets,
  $(n, \varepsilon)$-generating sets
  and $\covercardinality{\balls{d}{\varepsilon}^n}$
  can be found in \cite{patrao_caldas:entropia}.

  \begin{definition}
    \label{def:bowen_pressure}
    For a dynamical system $\dynamicalsystem{X}{T}$,
    a function $\function{f}{X}{\reals}$,
    $n = 1, 2, \dotsc$
    and $\varepsilon > 0$,
    define
    \begin{align*}
      \generatingpathpressure[n]{d}{T}{f}{\varepsilon}
      &=
      \inf
      \setsuchthat{\sum_{x \in E} \expbase^{f_n(x)}}{\text{$E$ is $(n,\varepsilon)$-generating}}
      \\
      \separatedpathpressure[n]{d}{T}{f}{\varepsilon}
      &=
      \sup
      \setsuchthat{\sum_{x \in E} \expbase^{f_n(x)}}{\text{$E$ is $(n,\varepsilon)$-separated}}
    \end{align*}
    and
    \begin{align*}
      \generatingpathpressure{d}{T}{f}{\varepsilon}
      &=
      \limsup_{n \rightarrow \infty}
      \frac{1}{n}
      \log \generatingpathpressure[n]{d}{T}{f}{\varepsilon}
      \\
      \separatedpathpressure{d}{T}{f}{\varepsilon}
      &=
      \limsup_{n \rightarrow \infty}
      \frac{1}{n}
      \log \separatedpathpressure[n]{d}{T}{f}{\varepsilon}.
    \end{align*}
    Since those last two are monotinic in $\varepsilon$,
    define
    \begin{align*}
      \generatingpressure{d}{T}{f}
      &=
      \lim_{\varepsilon \rightarrow 0}
      \generatingpathpressure{d}{T}{f}{\varepsilon}
      =
      \sup_{\varepsilon > 0}
      \generatingpathpressure{d}{T}{f}{\varepsilon}
      \\
      \separatedpressure{d}{T}{f}
      &=
      \lim_{\varepsilon \rightarrow 0}
      \separatedpathpressure{d}{T}{f}{\varepsilon}
      =
      \sup_{\varepsilon > 0}
      \separatedpathpressure{d}{T}{f}{\varepsilon}.
    \end{align*}
  \end{definition}

  We now state some very basic properties satisfied
  by the different kinds of topological pressure we have defined.
  First,
  let's relate them all.

  \begin{lemma}
    \label{lemma:topological_and_metric:easy_inequality}
    For a dynamical system $\dynamicalsystem{X}{T}$,
    any function $\function{f}{X}{\reals}$
    and any metric $d$ over $X$,
    \begin{equation*}
      \minustopologicalpressure{T}{f}
      \leq
      \plustopologicalpressure{T}{f}
      \leq
      \generatingpressure{d}{T}{f}
      \leq
      \separatedpressure{d}{T}{f}
      \leq
      \dtopologicalpressure{d}{T}{f}.
    \end{equation*}
  \end{lemma}

  \begin{proof}
    It is quite evident that
    $\minustopologicalpressure{T}{f} \leq \plustopologicalpressure{T}{f}$.
    The fact that
    $\generatingpressure{d}{T}{f} \leq \separatedpressure{d}{T}{f}$
    is a consequence of the fact that any $(n, \varepsilon)$-separated set
    is contained in a maximal one.
    And a maximal $(n,\varepsilon)$-separated set
    is in fact an $(n, \varepsilon)$-generating set.

    \proofitem
    {
      $\plustopologicalpressure{T}{f} \leq \generatingpressure{d}{T}{f}$
    }
    {
      Let $\family{A}$ be an admissible cover of $X$.
      Then,
      Lemma \ref{lemma:lebesgue_number} gives us
      $\varepsilon_0 > 0$ such that
      $\family{A} \prec \balls{d}{\varepsilon}$
      for any $\varepsilon \leq \varepsilon_0$.
      Let $E$ be an $(n, \varepsilon)$-generating set.
      In this case,
      \begin{equation*}
        \family{B}
        =
        \setsuchthat
        {
          \ball[d]{x}{\varepsilon}
          \cap \dotsb \cap
          T^{-(n-1)}\ball[d]{T^{n-1}x}{\varepsilon}
        }{x \in E}
        \setsuchthat{\ball[\iteratedmetric{d}{n}]{x}{\varepsilon}}{x \in E}
      \end{equation*}
      is a subcover of $\balls{d}{\varepsilon}^n$.
      Therefore,
      \begin{align*}
        \Qcoverpressure{n}{T}{f}{\family{A}}
        &\leq
        \Qcoverpressure{n}{T}{f}{\balls{d}{\varepsilon}}
        \\
        &\leq
        \sum_{B \in \family{B}}
        \inf \expbase^{f_n(B)}
        \\
        &\leq
        \sum_{x \in E}
        \expbase^{f_n(x)}.
      \end{align*}
      Taking the infimum for every $(n, \varepsilon)$-generating $E$,
      \begin{equation*}
        \Qcoverpressure{n}{T}{f}{\family{A}}
        \leq
        \generatingpathpressure[n]{d}{T}{f}{\varepsilon}.
      \end{equation*}
      Taking the logarithm, dividing by $n$,
      and taking the $\limsup$ for $n \rightarrow \infty$,
      we get that
      \begin{equation*}
        \plustcoverpressure{T}{f}{\family{A}}
        \leq
        \generatingpathpressure{d}{T}{f}{\varepsilon},
      \end{equation*}
      for every $\varepsilon < \varepsilon_0$.
      Therefore,
      by making $\varepsilon \rightarrow 0$,
      \begin{equation*}
        \plustcoverpressure{T}{f}{\family{A}}
        \leq
        \generatingpressure{d}{T}{f}.
      \end{equation*}
      Finally,
      since $\family{A}$ was an arbitrary admissible cover,
      \begin{equation*}
        \plustopologicalpressure{T}{f}
        \leq
        \generatingpressure{d}{T}{f}.
      \end{equation*}
    }

    Let us demonstrate the last inequality.
    \proofitem
    {
      $\separatedpressure{d}{T}{f} \leq \dtopologicalpressure{d}{T}{f}$
    }
    {
      Given $\varepsilon > 0$,
      let $E$ be any $(n, \varepsilon)$-separated set,
      and let $\family{B}$ be any subcover of $\balls{d}{\frac{\varepsilon}{2}}^n$.
      Then,
      for each $x \in E$,
      pick a $B_x \in \family{B}$,
      such that $x \in B_x$.
      Notice that,
      since $E$ is $(n, \varepsilon)$-separated,
      $x \neq y \Rightarrow B_x \neq B_y$.
      Therefore,
      \begin{align*}
        \sum_{x \in E} \expbase^{f_n(x)}
        &\leq
        \sum_{x \in E} \sup \expbase^{f_n(B_x)}
        \\
        &\leq
        \sum_{B \in \family{B}} \sup \expbase^{f_n(B)}.
      \end{align*}
      Taking the infimum for $\family{B} \subset \balls{d}{\frac{\varepsilon}{2}}$,
      and then the supremum for $(n, \varepsilon)$-separated $E$ gives
      \begin{equation*}
        \separatedpathpressure[n]{d}{T}{f}{\varepsilon}
        \leq
        \tPcoverpressure{T}{f}{\balls{d}{\frac{\varepsilon}{2}}^n}.
      \end{equation*}
      And the result follows by taking the logarithm,
      dividing by $n$,
      making $n \rightarrow \infty$,
      and then,
      taking the $\liminf$ for $\varepsilon \rightarrow 0$.
    }
  \end{proof}

  \begin{proposition}
    \label{proposition:topological_pressure:all_equal}
    For a dynamical system $\dynamicalsystem{X}{T}$,
    if $d$ is a one-point metric for $X$ and
    $\function{f}{X}{\reals}$ is one-point uniformly continuous,
    then
    \begin{equation*}
      \minustopologicalpressure{T}{f}
      =
      \plustopologicalpressure{T}{f}
      =
      \generatingpressure{d}{T}{f}
      =
      \separatedpressure{d}{T}{f}
      =
      \dtopologicalpressure{d}{T}{f}.
    \end{equation*}
  \end{proposition}

  \begin{proof}
    This is a direct consequence of Lemmas
    \ref{lemma:topological_and_metric:hard_inequality}
    and
    \ref{lemma:topological_and_metric:easy_inequality}.
  \end{proof}

  \begin{definition}[Topological Pressure]
    \label{def:topological_pressure}
    Let $\dynamicalsystem{X}{T}$ be a dynamical system
    that admits a metrizable one-point compactification.
    Suppose that $\function{f}{X}{\reals}$
    is one-point uniformly continuous.
    Then,
    the \emph{topological pressure}
    is the quantity in Proposition \ref{proposition:topological_pressure:all_equal},
    and is denoted by
    $\topologicalpressure{T}{f}$.
  \end{definition}

  When the space is compact
  and $\function{f}{X}{\reals}$ is continuous,
  it is a simple fact that
  \begin{equation*}
    \topologicalpressure{T^k}{f_k}
    =
    k \topologicalpressure{T}{f}.
  \end{equation*}
  For the non-compact case,
  only one inequality follows from a similar argument,
  and only for
  $\plustopologicalpressure{T}{f}$
  and
  $\minustopologicalpressure{T}{f}$.

  \begin{proposition}
    \label{proposition:sistema_iterado}
    Consider the dynamical system $\dynamicalsystem{X}{T}$,
    and a function $\function{f}{X}{\reals}$.
    Then,
    for any $k = 1, 2, \dotsc$,
    \begin{align*}
      \minustopologicalpressure{T^k}{f_k}
      &\leq
      k \minustopologicalpressure{T}{f}
      \\
      \plustopologicalpressure{T^k}{f_k}
      &\leq
      k \plustopologicalpressure{T}{f}.
    \end{align*}
  \end{proposition}

  \begin{proof}
    Let $\family{A}$ be an admissible cover of $X$.
    Notice that
    $(\family{A}_T^k)_{T^k}^n = \family{A}_T^{kn}$.
    And also,
    $(f_k)_{T^k,n} = f_{nk}$.
    So,
    \begin{equation*}
      \frac{1}{n}
      \log \Qcoverpressure{n}{T^k}{f_k}{\family{A}^k}
      =
      k\frac{1}{kn}
      \log \Qcoverpressure{kn}{T}{f}{\family{A}}.
    \end{equation*}
    Taking the $\liminf$ and $\limsup$ for $n \rightarrow \infty$,
    \begin{align*}
      \minustcoverpressure{T^k}{f_k}{\family{A}^k}
      &=
      k \minustcoverpressure{T}{f}{\family{A}}
      \\
      \plustcoverpressure{T^k}{f_k}{\family{A}^k}
      &=
      k \plustcoverpressure{T}{f}{\family{A}}.
    \end{align*}
    And since $\family{A} \prec \family{A}^k$,
    Lemma \ref{lemma:refinement:inequality} implies that
    \begin{align*}
      \minustcoverpressure{T^k}{f_k}{\family{A}}
      &\leq
      \minustcoverpressure{T^k}{f_k}{\family{A}^k}
      =
      k \minustcoverpressure{T}{f}{\family{A}}
      \\
      \plustcoverpressure{T^k}{f_k}{\family{A}}
      &\leq
      \plustcoverpressure{T^k}{f_k}{\family{A}^k}
      =
      k \plustcoverpressure{T}{f}{\family{A}}.
    \end{align*}
    Now, we just have to take the supremum for every admissible cover
    $\family{A}$ to reach the desired conclusion.
  \end{proof}

  Notice that the power of Proposition \ref{proposition:sistema_iterado}
  is quite limited even in the case where $\function{f}{X}{\reals}$
  is one-point uniformly continuous.
  In this case,
  eventhough $\topologicalpressure{T}{f}$ is defined,
  $\topologicalpressure{T^n}{f_n}$ might not be,
  because $f_n$ might not be one-point uniformly continuous.

  Let us finally mention a feature that is common
  to every concept of pressure we have defined so far.

  \begin{lemma}
    \label{lemma:pressure_plus_constant}
    If
    $d$ is any metric over $X$,
    $\function{f}{X}{\reals}$
    is any function
    and $c \in \reals$.
    Then,
    \begin{align*}
      \plustopologicalpressure{T}{f + c}
      &=
      \plustopologicalpressure{T}{f}
      + c
      \\
      \minustopologicalpressure{T}{f + c}
      &=
      \minustopologicalpressure{T}{f}
      + c
      \\
      \dtopologicalpressure{d}{T}{f + c}
      &=
      \dtopologicalpressure{d}{T}{f}
      + c
      \\
      \generatingpressure{d}{T}{f + c}
      &=
      \generatingpressure{d}{T}{f}
      + c
      \\
      \separatedpressure{d}{T}{f + c}
      &=
      \separatedpressure{d}{T}{f}
      + c.
    \end{align*}
    And if $\mu$ is a $T$-invariant probability measure
    and $f$ has a well defined integral,
    \begin{equation*}
      \measurepressure{\mu}{T}{f + c}
      =
      \measurepressure{\mu}{T}{f}
      + c.
    \end{equation*}
  \end{lemma}

  \begin{proof}
    For $\measurepressure{\mu}{T}{f + c}$,
    this is an obvious consequence of
    $\integral{(f + c)}{\mu} = \integral{f}{\mu} + c$.
    The other equalities are easy consequences of
    the exponential function properties.
  \end{proof}
}
}
  {\section{Variational Principle}
  \label{sec:principle}

  Inspired by what has been done for the compact case,
  we demonstrate a \emph{variational principle} for the pressure
  of a topological system $\function{T}{X}{X}$,
  where $X$ is not assumed to be compact
  but it is just assumed to have a one-point compactification $\onepointcompactification{X}$.
  This does not imply that $T$ can be itself extended
  to a topological dynamical system over $\onepointcompactification{X}$.

  We use the preparations made in Section
  \ref{section:preliminaries}
  in order to adapt Misiurewicz's demonstration
  of the variational principle.
  Misiurewicz's original article is \cite{misiurewicz}.
  We shall follow the more didatic presentation
  of the variational principle
  presented in \cite{viana:ergodic_theory}, Section $10.3$ and Section $10.4$.
  A similar presentation can also be found in
  \cite{walters}, Chapter $9$.

  We are concerned about the supremum of
  $\measurepressure{\mu}{T}{f}$
  over all $T$-invariant Radon probability measures
  for a given one-point uniformly continuous
  $\function{f}{X}{\reals}$.
  However,
  there might happen that no such a probability measure exists.
  In this case,
  we agree that
  \begin{equation*}
    \sup_{\mu}
    \measurepressure{\mu}{T}{f}
    =
    0.
  \end{equation*}
  According to Lemma \ref{lemma:entropia_na_compactificacao},
  this is the same as taking the supremum over all $T$-invariant
  Radon measures $\mu$ with $0 \leq \mu(X) \leq 1$.
  In this case, there is always an invariant measure.
  Namely, $\mu = 0$.

  \begin{theorem}
    \label{th:principio_variacional}
    Let $\dynamicalsystem{X}{T}$
    be a metrizable locally compact separable dynamical system,
    and let
    $\function{f}{X}{\reals} \in \vanishingatinfinityfunctions{X}$.
    Then,
    \begin{equation*}
      \sup_{\mu}
      \measurepressure{\mu}{T}{f}
      =
      \topologicalpressure{T}{f},
    \end{equation*}
    where the supremum is taken
    over all $T$-invariant Radon probability measures.
    If there is no $T$-invariant Radon probability measure,
    \begin{equation*}
      \topologicalpressure{T}{f}
      =
      0.
    \end{equation*}
  \end{theorem}

  Before the proof,
  let's extend Theorem \ref{th:principio_variacional} to
  one-point uniformly continuous functions.

  \begin{corollary}
    Let $\dynamicalsystem{X}{T}$
    be a metrizable locally compact separable dynamical system,
    and let
    $\function{f}{X}{\reals}$
    be one-point uniformly continuous.
    Then,
    \begin{equation*}
      \sup_{\mu}
      \measurepressure{\mu}{T}{f}
      =
      \topologicalpressure{T}{f},
    \end{equation*}
    where the supremum is taken
    over all $T$-invariant Radon probability measures.
    If there is no $T$-invariant Radon probability measure,
    \begin{equation*}
      \topologicalpressure{T}{f}
      =
      f(\infty).
    \end{equation*}
  \end{corollary}

  \begin{proof}
    Use the theorem with $f - f(\infty)$ in place of $f$.
    Then,
    use Lemma \ref{lemma:pressure_plus_constant}.
  \end{proof}

  The theorem will be demonstrated if we show that:
  \begin{enumerate}
    \item
      For any $T$-invariant Radon probability $\mu$,
      \begin{equation*}
        \measurepressure{\mu}{T}{f}
        \leq
        \topologicalpressure{T}{f}.
      \end{equation*}

    \item
      If we fix a one-point metric $d$,
      then,
      for any $\varepsilon > 0$,
      there is a $T$-invariant Radon measure $\mu$,
      with $0 \leq \mu(X) \leq 1$,
      such that
      \begin{equation*}
        \separatedpathpressure{d}{T}{f}{\varepsilon}
        \leq
        \measurepressure{\mu}{T}{f}.
      \end{equation*}
  \end{enumerate}
  These claims are the contents of the following two subsections.

  {\subsection{Topological Pressure is an Upper Bound}

  This subsection is devoted to the proof
  of the following proposition.
  The technique we present is a mix of what is done
  for Lemma 3.2 of \cite{patrao_caldas:entropia}
  and what is done in Section 10.4.1 in \cite{viana:ergodic_theory}.

  \begin{proposition}
    Let $\dynamicalsystem{X}{T}$ be a dynamical system
    such that $X$ has a metrizable one-point compactification,
    let $\function{f}{X}{\reals}$ be one-point uniformly continuous,
    and $\mu$ a $T$-invariant Radon probability measure.
    Then,
    \begin{equation*}
      \measurepressure{\mu}{T}{f}
      \leq
      \topologicalpressure{T}{f}.
    \end{equation*}
  \end{proposition}

  \begin{proof}
    Let $\mu$ be any $T$-invariant Radon probability measure.
    We shall show that
    for any $n = 1, 2, \dotsc$,
    \begin{equation}
      \label{eq:lemma:principio_variacional:desigualdade_facil:primeira_desigualdade}
      \measurepressure{\mu}{T^n}{f_n}
      \leq
      n \topologicalpressure{T}{f}
      +
      2
      +
      \log 2.
    \end{equation}
    And then,
    Proposition \ref{proposition:sistema_iterado:medida}
    implies that
    \begin{align*}
      \measurepressure{\mu}{T}{f}
      &=
      \frac{1}{n}
      \measurepressure{\mu}{T^n}{f_n}
      \\
      &\leq
      \topologicalpressure{T}{f}
      +
      \frac{2 + \log 2}{n}
      \rightarrow
      \topologicalpressure{T}{f}.
    \end{align*}
    And this will finish the demonstration.
    Notice that $f_n$ might not be one-point uniformly continuous,
    and therefore, we do not talk about
    $\topologicalpressure{T^n}{f_n}$.
    From now on,
    we fix $n$
    and attempt to show the validity of
    inequation \refeq{eq:lemma:principio_variacional:desigualdade_facil:primeira_desigualdade}.

    According to Proposition \ref{proposition:pressure_with_admissible_partitions},
    we have to show that
    given an admissible partition $\family{K}$,
    \begin{equation*}
      \tpartitionpressure{\mu}{T^n}{f_n}{\family{K}}
      \leq
      n \topologicalpressure{T}{f}
      +
      2
      +
      \log 2.
    \end{equation*}
    To that end,
    let $d$ be a one-point metric,
    it is enough if we prove that there is an $\varepsilon > 0$ such that
    \begin{equation}
      \label{eq:lemma:principio_variacional:desigualdade_facil:particao_e_cobertura}
      \tpartitionpressure{\mu}{T^n}{f_n}{\family{K}}
      \leq
      n \plustcoverpressure{T}{f}{\balls{d}{\varepsilon}}
      +
      2
      +
      \log 2.
    \end{equation}

    Let $\family{A}$ be the strongly admissible cover
    from Lemma \ref{lemma:cover_from_partition}.
    Using the \emph{Lebesgue Number} of Lemma \ref{lemma:lebesgue_number},
    fix $\varepsilon > 0$ such that
    \begin{equation*}
      \family{A}
      \prec
      \balls{d}{\varepsilon}.
    \end{equation*}
    Also,
    choose $\varepsilon$ small enough such that
    \begin{equation*}
      d(x,y)
      <
      2 \varepsilon
      \Rightarrow
      \abs{f(y) - f(x)}
      \leq
      \frac{1}{n}.
    \end{equation*}

    With $\varepsilon > 0$ properly choosen,
    we attempt at demonstrating the validity
    of inequality \refeq{eq:lemma:principio_variacional:desigualdade_facil:particao_e_cobertura}.
    Since we are working with $T$ and $T^n$ at the same time,
    let's agree that whenever the transformation is omitted,
    it is assumed to be $T$.

    \subproof
    {
      For any $m = 1, 2, \dotsc$,
      \begin{align*}
        \integral{f_n}{\mu}
        &+ 
        \frac{1}{m}
        \partitionentropy{\mu}{\family{K}_{T^n}^m}
        \leq \\ &\leq
        \frac{m+1}{m}
        +
        \log 2
        +
        \frac{n}{mn}
        \log \Qcoverpressure{mn}{T}{f}{\balls{d}{\varepsilon}}.
      \end{align*}
    }
    {
      Let
      $\family{B} \subset \balls{d}{\varepsilon}^{mn}$
      be any subcover.
      And notice that
      \begin{align*}
        \family{K}_{T^n}^m
        &\prec
        \family{K}^{mn}
        \\
        (f_n)_{T^n,m}
        &=
        f_{mn}.
      \end{align*}

      Given $C \in \family{K}_{T^n}^m$,
      let $x_C \in C$ be such that
      \begin{equation*}
        \sup f_{mn}(C)
        \leq
        f_{mn}(x_C) + 1.
      \end{equation*}
      Also,
      for each $C \in \family{K}_{T^n}^m$,
      choose $B_C \in \family{B}$ such that $x_C \in B_C$.
      Notice that for any $x \in B_C$
      and $j = 0, \dotsc, mn-1$,
      \begin{equation*}
        d(T^j x_C, T^j x)
        <
        2 \varepsilon.
      \end{equation*}
      Therefore,
      by the choice of $\varepsilon$,
      \begin{align*}
        \sup f_{mn}(C)
        &\leq
        f_{mn}(x_C)
        +
        1
        \\
        &\leq
        \inf f_{mn}(B_C)
        +
        \frac{mn}{n}
        +
        1
        \\
        &=
        \inf f_{mn}(B_C)
        +
        m
        +
        1.
      \end{align*}

      For each $B \in \family{B}$,
      let $c_B$ be the cardinality of
      \begin{equation*}
        \setsuchthat{C \in \family{K}_{T^n}^m}{B_C = B}.
      \end{equation*}
      Since $\family{A} \prec \balls{d}{\varepsilon} \prec \family{B}$,
      Lemma \ref{lemma:cover_from_partition:intersection_counting} implies that
      \begin{equation*}
        c_B
        \leq
        2^m.
      \end{equation*}

      Now,
      Lemma \ref{lemma:measure_pressure:upper_bound}
      with $T^n$ in place of $T$ and $f_n$ in place of $f$
      implies that
      \begin{align*}
        \integral{f_n}{\mu}
        &+ 
        \frac{1}{m}
        \partitionentropy{\mu}{\family{K}_{T^n}^m}
        \leq
        \frac{1}{m}
        \log \sum_{C \in \family{K}_{T^n}^m} \expbase^{\sup f_{mn}(C)}
        \\
        &\leq
        \frac{1}{m}
        \log
        \left(
          \expbase^{m+1} \sum_{C \in \family{K}_{T^n}^m} \expbase^{\inf f_{mn}(B_C)}
        \right)
        \\
        &=
        \frac{m + 1}{m}
        +
        \frac{1}{m}
        \log \sum_{C \in \family{K}_{T^n}^m} \expbase^{\inf f_{mn}(B_C)}
        \\
        &=
        \frac{m+1}{m}
        +
        \frac{1}{m}
        \log \sum_{B \in \family{B}} c_B\expbase^{\inf f_{mn}(B)}
        \\
        &\leq
        \frac{m+1}{m}
        +
        \frac{1}{m}
        \log
        \left(
          2^m \sum_{B \in \family{B}} \expbase^{\inf f_{mn}(B)}
        \right)
        \\
        &=
        \frac{m+1}{m}
        +
        \log 2
        +
        \frac{n}{mn}
        \log \sum_{B \in \family{B}} \expbase^{\inf f_{mn}(B)}.
      \end{align*}
      Taking the infimum for every subcover
      $\family{B} \subset \balls{d}{\varepsilon}^{mn}$,
      gives the Claim.
    }

    Now,
    use the Claim and take the $\limsup$ for $m \rightarrow \infty$
    \begin{align*}
      \integral{f_n}{\mu}
      +
      \tpartitionentropy{\mu}{\family{K}}{T^n}
      &\leq
      1 + \log 2 +
      \limsup_{m \rightarrow \infty}
      \frac{n}{mn}
      \log \Qcoverpressure{mn}{T}{f}{\balls{d}{\varepsilon}}
      \\
      &\leq
      1 + \log 2 +
      n \limsup_{k \rightarrow \infty}
      \frac{1}{k}
      \log \Qcoverpressure{k}{T}{f}{\balls{d}{\varepsilon}}
      \\
      &=
      1 + \log 2 +
      n \plustcoverpressure{T}{f}{\balls{d}{\varepsilon}}
      \\
      &\leq
      1 + \log 2 +
      n \plustopologicalpressure{T}{f}
      \\
      &=
      1 + \log 2 +
      n \topologicalpressure{T}{f},
    \end{align*}
    to get inequality \refeq{eq:lemma:principio_variacional:desigualdade_facil:particao_e_cobertura}
    and conclude the proof.
  \end{proof}
}
  {\subsection{Topological Pressure is a Lower Bound}

  This subsection is devoted to the proof
  of the following proposition,
  which is nothing more than
  a straight forward adaption of what is done
  in Subsection 10.4.2 of \cite{viana:ergodic_theory},
  using the same technique applied for Theorem 3.1 in \cite{patrao_caldas:entropia}.

  \begin{proposition}
    Let $\dynamicalsystem{X}{T}$ be a dynamical system
    such that $X$ admits a one-point compactification.
    Suppose
    $f \in \vanishingatinfinityfunctions{X}$.
    Then,
    for any $\varepsilon > 0$,
    there exists a $T$-invariant Radon measure $\mu$,
    with $0 \leq \mu(X) \leq 1$,
    such that
    \begin{equation*}
      \separatedpathpressure{d}{T}{f}{\varepsilon}
      \leq
      \measurepressure{\mu}{T}{f}.
    \end{equation*}
  \end{proposition}

  \begin{proof}
    Use Lemma \ref{lemma:compactificacao}
    to get a compact metrizable extension
    $\function{S}{Z}{Z}$ for $T$.
    According to Lemma \ref{lemma:entropia_na_compactificacao},
    the demonstration will be complete
    if we find a probability measure $\mu$ over $Z$
    which is $S$-invariant,
    and a partition $\family{C}$
    having a $C \in \family{C}$
    such that $\complementset{X} \subset C$,
    and such that
    \begin{equation}
      \label{proposition:hard_inequality:eq:separado_medida}
      \separatedpathpressure{d}{T}{f}{\varepsilon}
      \leq
      \tpartitionpressure{\mu}{S}{g}{\family{C}},
    \end{equation}
    where $g = f \circ \pi$,
    and
    $\function{\pi}{Z}{\onepointcompactification{X}}$ is the projection
    from Lemma \ref{lemma:compactificacao}.
    Notice that
    \begin{equation*}
      \integral{g}{\mu}
      =
      \integral{f}{\mu},
    \end{equation*}
    because $g|_{\complementset{X}} = f(\infty) = 0$.

    Let $d$ be a one-point metric for $X$,
    and $\widetilde{d}$ be the pseudometric over $Z$ induced by it.
    That is,
    considering $d$ as a metric over $\onepointcompactification{X}$,
    \begin{equation*}
      \widetilde{d}(x,y)
      =
      d(\pi(x), \pi(y)).
    \end{equation*}
    For each $n = 1, 2, \dotsc$,
    let $E_n \subset X$ be an $(n, \varepsilon)$-separated set such that
    \begin{equation*}
      \frac{1}{2}
      \separatedpathpressure[n]{d}{T}{f}{\varepsilon}
      \leq
      \sum_{x \in E_n}
      \expbase^{f_n(x)}.
    \end{equation*}
    Call the rightside quantity $A_n$.
    That is,
    \begin{equation*}
      \frac{1}{2}
      \separatedpathpressure[n]{d}{T}{f}{\varepsilon}
      \leq
      A_n.
    \end{equation*}

    {  Then,
  define over $Z$ the measure
  \begin{equation*}
    \sigma_n
    =
    \frac{1}{A_n}
    \sum_{x \in E_n}
    \expbase^{g_n(x)}
    \delta_x,
  \end{equation*}
  where $\delta_x$ is the Dirac measure with support in $x$.
  And notice that $\sigma_n$ is a probability measure.
  Also define
  \begin{equation*}
    \mu_n = \frac{1}{n} \sum_{j=0}^{n-1} \sigma_n \circ S^{-j}.
  \end{equation*}

  \subproof
  {
    There is a subsequence $n_k$
    and a Radon probability measure $\mu$
    such that $\mu_{n_k} \rightarrow \mu$,
    and such that
    \begin{equation*}
      \lim_{k \rightarrow \infty}
      \frac{1}{n_k}
      \log A_{n_k}
      =
      \limsup_{n \rightarrow \infty}
      \frac{1}{n}
      \log A_n.
    \end{equation*}
    Also,
    for any measurable $C \subset Z$
    with $\mu(\partial C) = 0$,
    \begin{equation*}
      \lim
      \mu_{n_k}(C)
      =
      \mu(C).
    \end{equation*}
  }
  {
    It is clear that there is a subsequence $n_k$ such that
    \begin{equation*}
      \lim_{k \rightarrow \infty}
      \frac{1}{n_k}
      \log A_{n_k}
      =
      \limsup_{n \rightarrow \infty}
      \frac{1}{n}
      \log A_n.
    \end{equation*}
    In the weak-$*$ topology,
    the set of Radon probability measures $\mu$ over $Z$
    is easily seen to be sequentially compact
    (Proposition 2.1.6 from \cite{viana:ergodic_theory}).

    The sequential compactness means that
    we can assume that $n_k$ is such that
    $\mu_{n_k}$ converges to some Radon probability $\mu$.
    The last assertion in our claim is a consequence of the
    \emph{Portmanteau Theorem},
    and can be found in \cite{billingsley:convergence},
    Theorem 2.1, item (v).
  }

  \subproof
  {
    The measure $\mu$ is $S$-invariant.
  }
  {
    It is clear that
    $\mu_{n_k} \circ S^{-1} \rightarrow \mu \circ S^{-1}$.
    In fact,
    for any continuous $\function{\phi}{Z}{\reals}$,
    $\phi \circ S$ is also continuous.
    Therefore,
    \begin{align*}
      \integral{\phi}{(\mu_{n_k} \circ S^{-1})}
      &=
      \integral{\phi \circ S}{\mu_{n_k}}
      \\
      &\rightarrow
      \integral{\phi \circ S}{\mu}
      \\
      &=
      \integral{\phi}{(\mu \circ S^{-1})}.
    \end{align*}
    On the other hand,
    \begin{align*}
      \abs{\integral{\phi}{(\mu_{n_k} - \mu_{n_k} \circ S^{-1})}}
      &=
      \abs{\integral{\phi \frac{1}{n_k}}{(\sigma_{n_k} - \sigma_{n_k} \circ S^{-n_k})}}
      \\
      &=
      \frac{1}{n_k}
      \abs{\integral{(\phi - \phi \circ S^{n_k})}{\sigma_{n_k}}}
      \\
      &\leq
      \frac{1}{n_k}
      \integral{\norm[\infty]{\phi - \phi \circ S^{n_k}}}{\sigma_{n_k}}
      \\
      &\leq
      \frac{1}{n_k}
      \integral{2 \norm[\infty]{\phi}}{\sigma_{n_k}}
      \\
      &=
      \frac{1}{n_k}
      2 \norm[\infty]{\phi}
      \rightarrow
      0.
    \end{align*}
    This implies that
    \begin{equation*}
      \mu
      =
      \lim \mu_{n_k}
      =
      \lim \mu_{n_k} \circ S^{-1}
      =
      \mu \circ S^{-1}.
    \end{equation*}
  }
}

    {  Now, we construct a suitable measurable partition $\family{Z}$,
  so that inequation \refeq{proposition:hard_inequality:eq:separado_medida} holds.
  To that end,
  we use the pseudometric $\widetilde{d}$.
  For each $z \in Z$,
  there exists a non null $\varepsilon_z < \frac{\varepsilon}{2}$
  such that the ball
  $B_z = \ball{\varepsilon_z}{z}$,
  centered at $z$ with radius $\varepsilon_z$,
  is such that $\mu(\partial B_z) = 0$.
  Such an $\varepsilon_z$ exists because
  since the border of the balls $\ball{\delta}{z}$ are all disjont,
  there is at most a countable number of
  reals $\delta < \frac{\varepsilon}{2}$ such that $\ball{\delta}{z}$
  has border with non null measure.
  Now,
  since $Z$ is compact and the balls are open,
  there is a finite number of such balls, $B_0, \dotsc, B_n$ covering $Z$.
  We can assume that $\set{B_0, \dotsc, B_n}$ has no proper subcover.
  Let
  \begin{equation*}
    Z_j
    =
    B_j
    \setminus
    \left(
      B_1
      \cup \dotsb \cup
      B_{j-1}
    \right).
  \end{equation*}
  Then,
  $\family{Z} = \set{Z_0, \dotsc, Z_k}$
  is a measurable partition.
  We can also assume that
  $\complementset{X} \subset B_0 = Z_0$,
  because in the pseudometric $\widetilde{d}$,
  $\complementset{X}$ has diameter equals to $0$.
  That is,
  $\family{Z}$ satisfies the condidtions of
  Lemma \ref{lemma:entropia_na_compactificacao}.

  Also,
  notice that each $C \in \family{Z}^n$
  is such that for any $x, y \in C$,
  \begin{equation*}
    \widetilde{d}(S^j x, S^j y)
    <
    \varepsilon
  \end{equation*}
  for all $j = 0, 1, \dotsc, n-1$.

  \subproof
  {
    For each
    $C \in \family{Z}^n$,
    $\mu(\partial C) = 0$.
  }
  {
    Notice that,
    since $S$ is continuous,
    the border operator $\partial$ possesses the following properties.
    \begin{enumerate}
      \item
        \label{it:border:property:complement}
        $\partial A = \partial \complementset{A}$.

      \item
        \label{it:border:property:intersection}
        $
          \partial (A_1 \cap \dotsb \cap A_k)
          \subset
          \partial A_1
          \cup \dotsb \cup
          \partial A_k
        $.

      \item
        \label{it:border:property:continuous}
        $
          \partial S^{-1}(A)
          \subset
          S^{-1}(\partial A)
        $.
    \end{enumerate}
    From items
    \refitem{it:border:property:complement}
    and
    \refitem{it:border:property:intersection},
    each $Z_j = B_j \cap \complementset{B_1} \cap \dotsb \cap \complementset{B_{j-1}}$
    in $\family{Z}$ has border with null measure.
    And from items
    \refitem{it:border:property:intersection}
    and
    \refitem{it:border:property:continuous},
    the same is true for the sets in $\family{Z}^n$.
  }
}

    Having constructed $\mu$ and $\family{C}$,
    it remains to show that
    inequation \refeq{proposition:hard_inequality:eq:separado_medida} holds.

    {  \subproof
  {
    $\integral{g_n}{\sigma_n} = n \integral{g}{\mu_n}$.
  }
  {
    In fact,
    \begin{align*}
      \integral{g}{\mu_n}
      &=
      \frac{1}{n}
      \sum_{j=0}^{n-1}
      \integral{g}{\sigma_n \circ S^{-j}}
      \\
      &=
      \frac{1}{n}
      \sum_{j=0}^{n-1}
      \integral{g \circ S^j}{\sigma_n}
      \\
      &=
      \frac{1}{n}
      \integral{\sum_{j=0}^{n-1} g \circ S^j}{\sigma_n}
      \\
      &=
      \frac{1}{n}
      \integral{g_n}{\sigma_n}.
    \end{align*}
  }

  \subproof
  {
    $\partitionentropy{\sigma_n}{\family{Z}^n}
    +
    n \integral{g}{\mu_n}
    =
    \log A_n$.
  }
  {
    Let $C \in \family{Z}^n$.
    Since each element of $\family{Z}$
    has diameter less then $\varepsilon$,
    we have that
    $C$ can contain at most one element $x \in E_n$.
    That is,
    $\sigma_n(C) = 0$ or
    $\sigma_n(C) = \frac{\expbase^{g_n(x)}}{A_n}$.
    Therefore,
    \begin{align*}
      \partitionentropy{\sigma_n}{\family{Z}^n}
      +
      n \integral{g}{\mu_n}
      &=
      \partitionentropy{\sigma_n}{\family{Z}^n}
      +
      \integral{g_n}{\sigma_n}
      \\
      &=
      \sum_{x \in E_n}
      \sigma_n(\set{x})
      \left(
        g_n(x)
        +
        \log \frac{1}{\sigma_n(\set{x})}
      \right)
      \\
      &=
      \sum_{x \in E_n}
      \frac{\expbase^{g_n(x)}}{A_n}
      \log \frac{\expbase^{g_n(x)}}{\expbase^{g_n(x)}/A_n}
      \\
      &=
      \sum_{x \in E_n}
      \frac{\expbase^{g_n(x)}}{A_n}
      \log A_n
      \\
      &=
      \log A_n.
    \end{align*}
  }

  Passing from $\sigma_n$ to $\mu_n$
  is the same procedure as in the compact case,
  as we shall detail right now.
  Notice that for any measurable finite partition $\family{D}$,
  Lemma 2.7 from \cite{patrao_caldas:entropia} implies that
  \begin{equation*}
    \sum_{j = 0}^{n-1}
    \frac{1}{n}
    \partitionentropy{\sigma_n \circ S^{-j}}{\family{D}}
    \leq
    \partitionentropy{\mu_n}{\family{D}}.
  \end{equation*}

  For
  $n, q \in \naturals$ with $1 < q < n$,
  take an integer $m$ such that $mq \geq n > m(q-1)$.
  Then, for every $j = 0, \dotsc, q-1$,
  \begin{align*}
    \family{Z}^n
    &\prec
    \family{Z}^j
    \vee
    S^{-j}
    \left(
      \family{Z}^{qm}
    \right)
    \\
    &=
    \family{Z}^j
    \vee
    S^{-j}(\family{Z}^q)
    \vee
    S^{-(j+q)}(\family{Z}^q)
    \vee \dotsb \vee
    S^{-(j+(m-1)q)}(\family{Z}^q).
  \end{align*}
  Therefore,
  using Lemma 2.6 from \cite{patrao_caldas:entropia},
  \begin{align*}
    \partitionentropy{\sigma_n}{\family{Z}^n}
    &\leq
    \partitionentropy{\sigma_n}{\family{Z}^j}
    +
    \partitionentropy{\sigma_n \circ S^{-(j + 0q)}}{\family{Z}^q}
    + \dotsb +
    \partitionentropy{\sigma_n \circ S^{-(j + (m-1)q)}}{\family{Z}^q}
    \\
    &\leq
    \partitionentropy{\sigma_n}{\family{Z}^q}
    +
    \partitionentropy{\sigma_n \circ S^{-(j + 0q)}}{\family{Z}^q}
    + \dotsb +
    \partitionentropy{\sigma_n \circ S^{-(j + (m-1)q)}}{\family{Z}^q}
    \\
    &\leq
    \log \cardinality{\family{Z}^q}
    +
    \partitionentropy{\sigma_n \circ S^{-(j + 0q)}}{\family{Z}^q}
    + \dotsb +
    \partitionentropy{\sigma_n \circ S^{-(j + (m-1)q)}}{\family{Z}^q}.
  \end{align*}
  Summing up in $j = 0, \dotsc, q-1$,
  \begin{align*}
    q \partitionentropy{\sigma_n}{\family{Z}^n}
    &\leq
    q \log \cardinality{\family{Z}^q}
    +
    \sum_{j=0}^{q-1}
    \sum_{a=0}^{m-1}
    \partitionentropy{\sigma_n \circ S^{-(j + aq)}}{\family{Z}^q}
    \\
    &=
    q \log \cardinality{\family{Z}^q}
    +
    \sum_{p=0}^{n-1}
    \partitionentropy{\sigma_n \circ S^{-p}}{\family{Z}^q}
    +
    \sum_{p=n}^{mq-1}
    \partitionentropy{\sigma_n \circ S^{-p}}{\family{Z}^q}
    \\
    &\leq
    2 q \log \cardinality{\family{Z}^q}
    +
    n
    \sum_{p=0}^{n-1}
    \frac{1}{n}
    \partitionentropy{\sigma_n \circ S^{-p}}{\family{Z}^q}
    \\
    &\leq
    2 q \log \cardinality{\family{Z}^q}
    +
    n \partitionentropy{\mu_n}{\family{Z}^q}.
  \end{align*}
  Since each element $C \in \family{Z}^q$ has border with null measure,
  \begin{equation*}
    \lim_{k \rightarrow \infty}
    \mu_{n_k}(C)
    =
    \mu(C).
  \end{equation*}
  An this implies that
  \begin{equation*}
    \partitionentropy{\mu_{n_k}}{\family{Z}^q}
    \rightarrow
    \partitionentropy{\mu}{\family{Z}^q}.
  \end{equation*}
  Therefore,
  \begin{align*}
    \limsup_{n \rightarrow \infty}
    \frac{1}{n}
    \log \separatedpathpressure[n]{d}{T}{f}{\varepsilon}
    &=
    \limsup_{n \rightarrow \infty}
    \frac{1}{n}
    \log \frac{1}{2} \separatedpathpressure[n]{d}{T}{f}{\varepsilon}
    \\
    &\leq
    \limsup_{n \rightarrow \infty}
    \frac{1}{n}
    \log A_n
    \\
    &=
    \lim_{k \rightarrow \infty}
    \frac{1}{n_k}
    \log A_{n_k}
    \\
    &=
    \lim_{k \rightarrow \infty}
    \frac{1}{n_k}
    \left(
      \partitionentropy{\sigma_{n_k}}{\family{Z}^{n_k}}
      +
      n_k \integral{g}{\mu_{n_k}}
    \right)
    \\
    &=
    \lim_{k \rightarrow \infty}
    \left(
      \frac{q \partitionentropy{\sigma_{n_k}}{\family{Z}^{n_k}}}{q n_k}
      +
      \integral{g}{\mu_{n_k}}
    \right)
    \\
    &\leq
    \lim_{k \rightarrow \infty}
    \left(
      \frac{1}{n_k}
      2 \log \cardinality{\family{Z}^q}
      +
      \frac{1}{q}
      \partitionentropy{\mu_{n_k}}{\family{Z}^q}
      +
      \integral{g}{\mu_{n_k}}
    \right)
    \\
    &=
    0
    +
    \frac{1}{q}
    \partitionentropy{\mu}{\family{Z}^q}
    +
    \integral{g}{\mu}
    \\
    &\xrightarrow{q \rightarrow \infty}
    \tpartitionentropy{\mu}{\family{Z}}{S}
    +
    \integral{g}{\mu}
    \\
    &=
    \tpartitionpressure{\mu}{S}{g}{\family{Z}}
    \\
    &\leq
    \measurepressure{\mu}{T}{f}.
  \end{align*}
  Where the last inequality is from
  Lemma \ref{lemma:entropia_na_compactificacao}.
}
  \end{proof}
}
}
}

  \bibliographystyle{aomalpha}
  \bibstyle{aomalpha}
  \bibliography{principio_variacional_para_pressao/principio_variacional_para_pressao}
\end{document}